\documentclass{amsart}
\usepackage{cases}
\usepackage{amsmath}
\usepackage{amssymb}
\usepackage{amsfonts}
\usepackage{graphicx}
\newtheorem{theorem}{Theorem}[section]

\newtheorem{claim}[theorem]{Claim}

\newtheorem{corollary}[theorem]{Corollary}

\newtheorem{lemma}[theorem]{Lemma}

\newtheorem*{problem*}{Problem 1}

\newtheorem{remark}[theorem]{Remark}

\newtheorem{thm}{Theorem}

\begin{document}
\title[Pinching estimates in higher dimensions]
{Pinching estimates for solutions of the linearized Ricci flow system\\
in higher dimensions}
\author{Jia-Yong Wu and Jian-Biao Chen}

\address{Department of Mathematics, Shanghai Maritime University,
1550 Haigang Avenue, Shanghai 201306, P. R. China}

\email{jywu81@yahoo.com}

\email{chenjb@shmtu.edu.cn}

\thanks{The first author is partially supported by the NSFC (11101267, 11271132)
and the Innovation Program of Shanghai Municipal Education Commission
(13YZ087). The second author is partially supported by the NSFC (51109128)}


\date{\today}

\dedicatory{}

\begin{abstract}
We prove pinching estimates for solutions of the linearized Ricci flow
system on a closed manifold of dimension $n\geq 4$ with positive
scalar curvature and vanishing Weyl tensor. If the vanishing
Weyl tensor condition is removed, we only give a rough pinching
estimate controlled by some blow-up function in a short time. These results
extend the $3$-dimensional case due to Anderson and Chow (2005) \cite{[AC]}.
\end{abstract}
\maketitle

\section{Introduction and main results}
Given an $n$-dimensional closed Riemannian manifold $M^n$, a
smooth family of Riemannian metrics $g(t)$, $t\in [0, T)$, is
said to be evolving under the \emph{Ricci flow} if
\begin{equation}\label{Rf}
\frac{\partial}{\partial t}g_{ij}=-2R_{ij},
\end{equation}
where $R_{ij}$ is the Ricci curvature of the metric $g(t)$.
This geometric flow was introduced by R. Hamilton \cite{[Ham1]}.
In \cite{[Ham1]}, he proved that for any an initial metric $g_0$,
there is a unique solution to this flow over some short time
interval, with $g(0)=g_0$. This proof was later simplified by
D. DeTurck \cite{[De]} by linearizing the modified Ricci flow,
which leads to the following \emph{Lichnerowicz Laplacian heat
equation}
\begin{equation}\label{LLh}
\frac{\partial}{\partial t}h_{ij}=(\Delta_Lh)_{ij}
:=\Delta h_{ij}+2R_{ikjl}h_{kl}-R_{ik}h_{kj}-R_{jk}h_{ki},
\end{equation}
for a symmetric $2$-tensor $h$, where $R_{ikjl}$ denotes the
Riemannian curvature of the metric $g(t)$ moving under the
Ricci flow. Nowadays the Ricci flow \eqref{Rf} coupled with
equation \eqref{LLh} is often called the \emph{linearized
Ricci flow system}, which is useful for studying the singularities
of the Ricci flow. If we let $h_{ij}=R_{ij}$, then equation
\eqref{LLh} is exactly the evolution of the Ricci curvature
under the Ricci flow.

On the subject of differential Harnack inequalities for the
linearized Ricci flow system, there have been many important
contributions. In \cite{[ChHam]}, B. Chow and R. Hamilton
proved a linear trace differential Harnack inequality for
the coupled system \eqref{Rf}-\eqref{LLh} by adapting
similar techniques of Hamilton's trace inequality for
the Ricci flow \cite{[Ham2]} and Li-Yau's seminal inequality
for the heat equation \cite{[Li-Yau]}. Meanwhile in
\cite{[ChCh]}, B. Chow and S.-C. Chu gave a geometric
interpretation in terms of this linear trace differential
Harnack inequality. In \cite{[GIK]}, C. Guenther,
J. Isenberg and D. Knopf studied the linearized Ricci flow
system at flat solutions from the point of view of maximal
regularity theory. Besides, the related Li-Yau-Hamilton type
differential Harnack inequalities for the linearized Ricci flow
system were also appeared in \cite {[Chow]}, \cite{[ChKn]} and
\cite{[WuZh]}.

In another direction, in \cite{[AC]}, G. Anderson and B. Chow
considered the quantity $|h|^2/R^2$ for solutions to the linearized
Ricci flow system \eqref{Rf}-\eqref{LLh} on closed $3$-manifolds
with positive scalar curvature. They proved an interesting pinching
estimate for solutions of the linearized Ricci flow system on a closed
$3$-manifold:
\begin{thm}[see Anderson and Chow \cite{[AC]}]\label{ACthm}
Let $(M^3,g(t))$ be a solution to the Ricci flow on a closed $3$-manifold
on a time interval $[0,T)$ with $T<\infty$ and let $\rho\in[0,\infty)$ be
such that $R_{min}(0)>-\rho$. If the pair $(g,h)$ is any solution to the
linearized Ricci flow system \eqref{Rf}-\eqref{LLh}, then there exists a
constant $C<\infty$ depending only on $g(0)$, $h(0)$, $\rho$ and $T$
such that
\[
\frac{|h|}{R+\rho}\leq C
\]
on $M\times[0,T)$. Furthermore, when $\rho=0$, $C$ is independent of $T$.
\end{thm}
The inspiration for the above pinching estimate partly comes from the works of
R. Hamilton (see Section 10 of \cite{[Ham1]} and Section 24 of \cite{[Ham3]})
and M. Gursky \cite{[Gu]}. We refer the reader to the excellent introduction
of \cite{[AC]} for nice explanations on this subject. In Theorem \ref{ACthm},
if we take $h=Rc$ and $\rho=0$, we immediately obtain Hamilton's Ricci
pinching estimate:
\begin{thm}[see Hamilton \cite{[Ham1]}]\label{Hamth}
If $(M^3,g(t))$, $t\in[0,T)$, $T<\infty$ is a solution to the Ricci flow
\eqref{Rf} on a closed $3$-manifold with positive scalar curvature,
then there exists a constant $C<\infty$ depending only on $g(0)$
such that
\[
\frac{|Rc|}{R}\leq C
\]
on $M\times[0,T)$.
\end{thm}

Recently, S. Brendle \cite{[Br1]} has successfully applied Theorem
\ref{ACthm} (see Proposition 6.1 in \cite{[Br1]}) to give a complete
proof of the uniqueness of the $3$-dimensional Bryant soliton
which is $k$-noncollapsed and non-flat. This resolves a well-known
question mentioned in Perelman's paper \cite{[Per]}. Later in
\cite{[Br2]}, Brendle generalized his unique result to the higher
dimensional setting. That is, he obtained a similar uniqueness
theorem for the steady gradient Ricci soliton of dimension
$n\geq 4$ which has positive sectional curvature and is
asymptotically cylindrical. In the proof of the higher dimensions,
Brendle adapted the arguments in \cite{[Br1]} but needed some
extra technique work. For example, due to the lack of an analogue of
Theorem \ref{ACthm} in higher dimensions, he employed a new
argument to prove an estimate related the Lichnerowicz-type
equation (see Proposition 4.2 in \cite{[Br2]}), which is quite different
from the $3$-dimensional case.

In this paper, we mainly generalize Anderson-Chow's pinching estimate
to the higher dimensions under some curvature assumptions. Our
results could be viewed as a partial answer to a question implied in
the paper \cite{[Br2]}, which may be expected to study the singularities
of the Ricci flow. The results of this paper can be divided into two
parts.

On one hand, we will establish an Anderson-Chow's pinching estimate
for solutions of the linearized Ricci flow system on a closed manifold
of dimension $n\geq 4$ as long as scalar curvature preserves positive
and the Weyl tensor remains identically zero under the linearized
Ricci flow system. The Weyl tensor is an important geometric quantity
in understanding higher dimensional Riemannian manifolds. This tensor
is known to depend on the conformal structure, so that if
$\tilde{g}_{ij}=\phi g_{ij}$, then $\tilde{W}_{ijkl}=\phi W_{ijkl}$.
If $n\leq 3$, the Weyl tensor vanishes; if $n\geq 4$, in local coordinate,
the Weyl tensor can be written as
\begin{equation}
\begin{aligned}\label{Weyl}
W_{ikjl}&=R_{ikjl}-\frac{1}{n-2}(g_{ij}R_{kl}+g_{kl}R_{ij}-g_{il}R_{jk}-g_{jk}R_{il})\\
&\quad+\frac{R}{(n-1)(n-2)}(g_{ij}g_{kl}-g_{il}g_{jk}).
\end{aligned}
\end{equation}

There exist some examples of the Ricci flow with vanishing Weyl tensors
for all time. One distinct example is the 1-parameter family of conformally
equivalent metrics
\[
g(t)=\left(1-2(n-1)t\right)g_{S^n}
\]
on the unit round sphere $\mathbb{S}^n$ with the standard metric $g_{S^n}$.
This is an ancient solution of the Ricci flow with vanishing Weyl tensors for
$-\infty<t<\frac{1}{2(n-1)}$. Another interesting example is the time-parameter
family of metrics
\[
g(t)=\left(1-2(n-2)t\right)g_{S^{n-1}}+dr^2
\]
on the round cylinder $\mathbb{S}^{n-1}\times \mathbb{R}$
with the metric $g_0=g_{S^{n-1}}+dr^2$, $n\geq 3$. It is easy to check
that it is an ancient solution of the Ricci flow
with vanishing Weyl tensors for $-\infty<t<\frac{1}{2(n-2)}$.

Now we give the first main result, which describes a precise pointwise measure of the size
of $h$ relative to the scalar curvature in higher dimensions when the Weyl tensor remains
identically zero under the Ricci flow.
\begin{theorem}\label{Main}
Let $(M^n,g(t),h(t))$, $t\in[0,T)$, $T<\infty$ be a solution to the
linearized Ricci flow system \eqref{Rf}-\eqref{LLh} on a closed
$n$-dimensional ($n\geq 4$) manifold with positive scalar curvature.
Assume that the Weyl tensor remains identically zero under the
linearized Ricci flow system. Then there exists a constant
$c_0<\infty$ depending only on $g(0)$ and $h(0)$ such that
\[
\frac{|h|}{R}(x,t)\leq c_0
\]
on $M\times[0,T)$.
\end{theorem}
The proof of Theorem \ref{Main} essentially follows the arguments of \cite{[AC]}.
The main difference is more complicated curvature terms due to the higher
dimensions. We reduce this difficulty to a matrix problem, which is discussed in
Section \ref{NonP}. The assumption of the Weyl tensor in Theorem
\ref{Main} seems to be necessary because the Weyl tensor remains
identically zero in Theorem A. It is interesting to ask if there exists
a similar pinching estimate under different curvature assumptions.

\begin{remark}
In Theorem \ref{Main}, we assume that the scalar curvature is positive; however
in Theorem A, the scalar curvature only needs to be a lower bound. This is
because Anderson and Chow use a special property of curvature estimate in
dimension 3 due to Hamilton-Ivey estimate (see Theorem 24.4 in \cite{[Ham3]}
or \cite{[Ivey]}), whereas this property does not hold in higher dimensions.
Therefore the scalar curvature condition in Theorem \ref{Main} is a
little stronger than the case of Theorem A.
\end{remark}

\begin{remark}
If we let $h=Rc$, then Theorem \ref{Main} is just as a
special case of Knopf's Ricci curvature pinching
estimate \cite{[Kno]} (see also \cite{[Caox]}).
\end{remark}
\begin{remark}
Recently, X.-D. Cao and H. Tran \cite{[CaoxTr]} observed
that the Ricci flow solution with nonnegative isotropic
curvature implies a Riemannian curvature pinching result:
\[
\frac{|Rm|}{R}(x,t)\leq c
\]
for some constant $c=c(n,g(0))<\infty$.
\end{remark}

\vspace{0.5em}

On the other hand, if the vanishing Weyl tensor condition is removed in
Theorem \ref{Main}, we can give another Anderson-Chow's type pinching estimate
for solutions of the linearized Ricci flow system on a closed manifold.
Though this pinching estimate is not uniformed by a constant, it could
be controlled by some blow-up function in a short time.
\begin{theorem}\label{Main2}
Let $(M^n,g(t),h(t))$, $t\in[0,T)$, $T<\infty$ be a solution to the
linearized Ricci flow system \eqref{Rf}-\eqref{LLh} on a closed
$n$-dimensional ($n\geq 4$) manifold with positive scalar curvature,
and let $K:=\max_{t=0}|Rm|$. Then there exist finite constants $c_0:=c_0(g(0),h(0))$
and $c:=c(n)$ such that
\[
\frac{|h|^2}{R^2}(t)\leq c_0\cdot(1-8Kt)^{-c}
\]
on $M\times[0,T')$, where $T':=\min\{T, \frac{1}{8K}\}$.
\end{theorem}

The rest of this paper organized as follows. In Section \ref{BaRe},
following the arguments of \cite{[AC]}, we will prove Theorem \ref{Main}
by the straightforward computation and the usage of parabolic maximum
principle. The main difference is that we need to estimate the
nonnegativity of some complicated terms in the evolution equation due
to the higher dimensions. In Section \ref{BaRe2}, we first give
an estimate on the norm of Riemannian curvature under the Ricci flow
on a closed manifold. Then we apply this estimate and parabolic
maximum principles to prove Theorem \ref{Main2}. In Section \ref{NonP},
we will use the basic matrix theory to justify the nonnegativity
of the above mentioned terms appeared in the proof of
Theorem \ref{Main}.

\vspace{0.5em}

\textbf{Acknowledgement}
The author would like to thank the referee for many valuable
comments and suggestions on an earlier version of this paper.


\section{Proof of Theorem \ref{Main}}\label{BaRe}
Along the Ricci flow \eqref{Rf}, we have
\begin{equation}\label{evolubian}
\frac{\partial }{\partial t}R=\Delta R+2|Rc|^2.
\end{equation}
Combining equations \eqref{Rf}, \eqref{LLh}, \eqref{evolubian},
and the key estimate (see Claim \ref{Klem}), we complete
the proof of Theorem \ref{Main}.
\begin{proof}[Proof of Theorem \ref{Main}]
The proof involves a direct computation and the parabolic
maximum principle. Here we can borrow Anderson-Chow's
computation to simplify our calculation. From
the evolution equation (16) in \cite{[AC]}, we have
\begin{equation}
\begin{aligned}\label{evolu}
\frac{\partial}{\partial t}\left(\frac{|h|^2}{R^2}\right)
&=\Delta\left(\frac{|h|^2}{R^2}\right)+\frac{2}{R}\nabla R\cdot\nabla\left(\frac{|h|^2}{R^2}\right)
-\frac{2}{R^4}\left|R\nabla_ih_{jk}-\nabla_iRh_{jk}\right|^2\\
&\quad+\frac{4}{R^2}R_{ikjl}h_{ij}h_{kl}-\frac{4}{R^3}|h|^2|Rc|^2.
\end{aligned}
\end{equation}
The above formula holds for all dimensions. Since here $n\geq 4$
and the Weyl tensor $W=0$ by assumption, from \eqref{Weyl}
we have
\begin{equation*}
\begin{aligned}
R_{ikjl}&=\frac{1}{n-2}(g_{ij}R_{kl}+g_{kl}R_{ij}-g_{il}R_{jk}-g_{jk}R_{il})\\
&\quad-\frac{R}{(n-1)(n-2)}(g_{ij}g_{kl}-g_{il}g_{jk}).
\end{aligned}
\end{equation*}
Then substituting this into \eqref{evolu} yields
\begin{equation}
\begin{aligned}\label{evolu2}
\frac{\partial}{\partial t}\left(\frac{|h|^2}{R^2}\right)
&=\Delta\left(\frac{|h|^2}{R^2}\right)
+\frac{2}{R}\nabla R\cdot\nabla\left(\frac{|h|^2}{R^2}\right)\\
&\quad-\frac{2}{R^4}\left|R\nabla_ih_{jk}-\nabla_iRh_{jk}\right|^2
-\frac{4}{R^3}P,
\end{aligned}
\end{equation}
where
\begin{equation*}
\begin{aligned}
P:=&|h|^2|Rc|^2-\frac{Rh_{ij}h_{kl}}{n-2}(g_{ij}R_{kl}+g_{kl}R_{ij}-g_{il}R_{jk}-g_{jk}R_{il})\\
&+\frac{R^2h_{ij}h_{kl}}{(n-1)(n-2)}(g_{ij}g_{kl}-g_{il}g_{jk})\\
=&|h|^2|Rc|^2-\frac{2R}{n-2}(HR_{ij}h_{ij}-R_{jk}h_{ji}h_{ik})+\frac{R^2(H^2-|h|^2)}{(n-1)(n-2)}
\end{aligned}
\end{equation*}
and where $H:=g^{ij}h_{ij}$.

Now we need to deal with the troublesome term $P$. Fortunately,
we observe that $P$ is always nonnegative
without any assumption, which shall be confirmed
in the Section \ref{NonP}
(see Corollary \ref{cor}).
\begin{claim}\label{Klem}
If $n\geq 4$, for any metric $g$ and symmetric $2$-tensor $h$, we have
\[
|h|^2|Rc|^2-\frac{2R}{n-2}(HR_{ij}h_{ij}-R_{jk}h_{ji}h_{ik})+\frac{R^2(H^2-|h|^2)}{(n-1)(n-2)}\geq 0.
\]
\end{claim}

Proceeding our proof, by Claim \ref{Klem}
and $R>0$, from \eqref{evolu2} we derive that
\[
\frac{\partial}{\partial t}\left(\frac{|h|^2}{R^2}\right)
\leq\Delta\left(\frac{|h|^2}{R^2}\right)
+\frac{2}{R}\nabla R\cdot\nabla\left(\frac{|h|^2}{R^2}\right).
\]
Finally, applying the parabolic maximum principle to the above equation
yields
\[
\frac{|h|^2}{R^2}(x,t)\leq c_0
\]
on $M\times[0,T)$, where $c_0:=\max_{t=0}|h|^2/R^2$.
\end{proof}
\begin{remark}
The theorem is also true on complete noncompact Riemannian manifolds
when we can apply the maximum principle.
\end{remark}

\section{Proof of Theorem \ref{Main2}}\label{BaRe2}
In this section, we will prove Theorem \ref{Main2}. First,
we give a Riemannian curvature estimate under the Ricci flow
on a closed manifold, which will be useful in the proof of
Theorem \ref{Main2}.
\begin{lemma}\label{Riemcon}
Let $(M^n,g(t))$ be a solution to the Ricci flow on a closed $n$-dimensional
manifold on a time interval $[0,T)$ with $T<\infty$. Then
\begin{equation}\label{Riem}
|Rm(x,t)|\leq \frac{K}{1-8Kt}
\end{equation}
on $M\times[0,T')$, where $T':=\min\{T, \frac{1}{8K}\}$ and
$K:=\max_{t=0}|Rm|$.
\end{lemma}

\begin{proof}
Under the Ricci flow,
\[
\frac{\partial}{\partial t}|Rm|^2
\leq\Delta|Rm|^2-2|\nabla Rm|^2+16|Rm|^3.
\]
By the maximum principle, we have
\[
|Rm(g(t))|\leq \frac{K}{1-8Kt}
\]
on $M\times[0,T')$, where $T':=\min\{T, \frac{1}{8K}\}$ and
$K:=\max_{t=0}|Rm|$.
\end{proof}

Now we use Lemma \ref{Riemcon} to finish the proof of Theorem
\ref{Main2}.
\begin{proof}[Proof of Theorem \ref{Main2}]
As before, we still have the evolution equation \eqref{evolu}.
Since $(M^n,g(t))$, $t\in[0,T)$, $T<\infty$ has a solution of the
Ricci flow \eqref{Rf} on a closed $n$-dimensional ($n\geq 4$)
manifold, by Lemma \ref{Riemcon}, we have
\begin{equation}\label{Rmest}
|Rm(x,t)|\leq\frac{K}{1-8Kt}
\end{equation}
on $M\times[0,T')$, where $T':=\min\{T, \frac{1}{8K}\}$ and
$K:=\max_{t=0}|Rm|$. Substituting \eqref{Rmest}
into \eqref{evolu} and noticing $R>0$ by our assumption, we have
\[
\frac{\partial}{\partial t}\left(\frac{|h|^2}{R^2}\right)
\leq\Delta\left(\frac{|h|^2}{R^2}\right)+\frac{2}{R}\nabla R\cdot\nabla\left(\frac{|h|^2}{R^2}\right)
+\frac{C(n)K}{1-8Kt}\cdot\frac{|h|^2}{R^2}.
\]
Applying the parabolic maximum principle to the above
evolution equation (for example, see Proposition 4.3
in \cite{[ChKn2]}), we conclude that
\[
\frac{|h|^2}{R^2}(t)\leq c_0\cdot(1-8Kt)^{-c}
\]
on $M\times[0,T')$, where $c_0:=\max_{t=0}|h|^2/R^2$ and $c:=c(n)$.
This finishes the proof of the theorem.
\end{proof}


\section{Nonnegativity of a degree $4$ homogeneous polynomial\\
 in $2n$ variables}\label{NonP}
In this section, we will prove Claim \ref{Klem} in Section \ref{BaRe}.
Our proof seems to be different from \cite{[AC]}. Adapting the
Anderson-Chow's symbols in \cite{[AC]}, since $h$ is a symmetric tensor,
we may assume $h$ is diagonal. Let $h_1$, $h_1$,$\ldots$, $h_n$ denote the
eigenvalues of $h$ and let $r_1=R_{11}$, $r_2=R_{22}$, $\ldots$, $r_n=R_{nn}$
denote the diagonal entries of $R_{ij}$. Then
\begin{equation}
\begin{aligned}\label{PQ}
P&=|Rc|^2|h|^2-\frac{2R}{n-2}(HR_{ij}h_{ij}-R_{jk}h_{ji}h_{ik})+\frac{R^2(H^2-|h|^2)}{(n-1)(n-2)}\\
&\geq Q:=\sum^n_{i=1}r^2_i\cdot\sum^n_{i=1}h^2_i\\
&\quad+\frac{2}{n-2}\sum^n_{i=1}r_i\cdot\left(-\sum^n_{i=1}h_i\sum^n_{i=1}r_ih_i+\sum^n_{i=1}r_ih^2_i\right)\\
&\quad+\frac{1}{(n-1)(n-2)}\left(\sum^n_{i=1}r_i\right)^2\cdot
\left[\left(\sum^n_{i=1}h_i\right)^2-\sum^n_{i=1}h^2_i\right],
\end{aligned}
\end{equation}
where we used the fact $|Rc|^2\geq\sum^n_{i=1}r^2_i$.

Now we rewrite $Q$ as a bilinear form in $\gamma=(r_1,r_2,\ldots,r_n)^T$ with
coefficients in $h=(h_1,h_2,\ldots,h_n)^T$:
\[
Q=\gamma^T(s_2I_n+\alpha_0\beta^T+\beta\alpha^T_0)\gamma,
\]
where $I_n$ is the identity matrix of order $n$,
\begin{equation}\label{def1}
\alpha_0:=(1,1,\ldots,1)^T\in \mathbb{R}^n,
\quad\quad
\alpha_k:=(h^k_1,h^k_2,\ldots,h^k_n)^T\in \mathbb{R}^n
\end{equation}
and
\begin{equation}\label{def2}
s_0=n,\quad\quad s_k:=\sum^n_{i=1}h^k_i
\end{equation}
for $k=1,2,3,4$; and where
\begin{equation}\label{def3}
\beta:=\frac{(s^2_1-s_2)\alpha_0}{2(n-1)(n-2)}-\frac{s_1\alpha_1}{n-2}+\frac{\alpha_2}{n-2}.
\end{equation}

In the rest of this section, we shall prove $Q\geq 0$. To achieve
it, we begin
with two technical lemmas.
\begin{lemma}\label{tech1}
If the matrix $\Lambda\in \mathbb{R}^{n\times n}$ is nonsingular, then for any
column vector $\xi_i,\eta_i\in\mathbb{R}^n$, $(i=1,2)$
\[
\det\left(\Lambda+\xi_1\eta_1^T+\xi_2\eta_2^T\right)
=\det(\Lambda)\cdot\det\left(I_2+(\eta_1,\eta_2)^T\Lambda^{-1}(\xi_1,\xi_2)\right).
\]
\end{lemma}

\begin{proof}
One can easily check that
\begin{align*}
\begin{bmatrix}
I_{2}    & ~~~(\eta_1,\eta_2)^T \\
0        & ~~~\Lambda+\xi_1\eta_1^T+\xi_2\eta_2^T
\end{bmatrix}
=&\begin{bmatrix}
I_{2}    & ~~~0 \\
(\xi_1,\xi_2)  & ~~~I_n
\end{bmatrix}
\times
\begin{bmatrix}
I_{2}    & ~~~(\eta_1,\eta_2)^T\Lambda^{-1}\\
0        & ~~~I_n
\end{bmatrix}\\
&\times\begin{bmatrix}
I_2+(\eta_1,\eta_2)^T\Lambda^{-1}(\xi_1,\xi_2)  &~~~ 0\\
-(\xi_1,\xi_2)  & ~~~\Lambda
\end{bmatrix}.
\end{align*}
Then the conclusion can be immediately obtained by the above
identity.
\end{proof}

\begin{lemma}\label{tech2}
Let $\alpha_0$, $\alpha_k$, $s_k$ and $\beta$ be defined by \eqref{def1},
\eqref{def2} and \eqref{def3}. Then
\begin{equation}\label{buden1}
s_2+\alpha_0^T\beta\geq0
\end{equation}
and
\begin{equation}\label{buden2}
(s_2+\alpha_0^T\beta)^2-n\beta^T\beta\geq0.
\end{equation}
\end{lemma}

\begin{proof}
We first check \eqref{buden1}. Clearly,
$\alpha^T_i\alpha_j=s_{i+j}$ for all $i,j=0,1,2$.
Using this, by the definition $\beta$, direct computation yields
\begin{equation*}
\begin{aligned}
s_2+\alpha_0^T\beta&=s_2+\frac{n(s_1^2-s_2)}{2(n-1)(n-2)}-\frac{s_1^2}{n-2}+\frac{s_2}{n-2}\\
&=\left(1+\frac{1}{2(n-1)}\right)s_2-\frac{s_1^2}{2(n-1)}\\
&=\alpha_1^T\left[\left(1+\frac{1}{2(n-1)}\right)I_n-\frac{\alpha_0\alpha_0^T}{2(n-1)}\right]\alpha_1.
\end{aligned}
\end{equation*}
The fact that $1>\sum^{n-1}_{i=1}\frac{1}{2(n-1)}=\frac 12$ tells us that the real
matrix
\[
A:=\left(1+\frac{1}{2(n-1)}\right)I_n-\frac{\alpha_0\alpha_0^T}{2(n-1)}
\]
is strictly diagonally dominant. Moreover all main diagonal entries of
$A$ are positive. Hence, by Theorem 6.1.10 of \cite{[HoJo]}, we conclude that
$A$ is positive definite and
\[
s_2+\alpha_0^T\beta=\alpha_1^TA\alpha_1\geq 0.
\]

Then, we shall prove \eqref{buden2}. According to the definitions
$\alpha_k$ and $s_k$, using the relation $\alpha^T_i\alpha_j=s_{i+j}$,
we expand the expression
$(s_2+\alpha_0^T\beta)^2-n\beta^T\beta$ as
\begin{equation*}
\begin{aligned}
&(s_2+\alpha_0^T\beta)^2-n\beta^T\beta
=\left[s_2+\frac{n(s_1^2-s_2)}{2(n-1)(n-2)}-\frac{s^2_1}{n-1}+\frac{s_2}{n-1}\right]^2\\
&-n\left[\frac{n(s_1^2-s_2)^2}{4(n-1)^2(n-2)^2}+\frac{s_1^2s_2+s_4-2s_1s_3}{(n-2)^2}
+\frac{(s_1^2-s_2)s_2-(s_1^2-s_2)s_1^2}{(n-1)(n-2)^2}\right].
\end{aligned}
\end{equation*}
From above, we easily judge that
\[
f(h_1,h_2,\ldots,h_n):=(s_2+\alpha_0^T\beta)^2-n\beta^T\beta
\]
is a real homogeneous symmetric polynomial of degree $4$
with respect to $n$ variables: $h_1,h_2,\ldots,h_n$. Now
we recall an interesting fact due to V. Timofte
(see Corollary 5.6 in \cite{[Tim]}).
\begin{thm}[see Timofte \cite{[Tim]}]\label{propTi}
If $f$ is a real homogeneous symmetric polynomial of degree $4$
on $\mathbb{R}^n$, then
\[
f(h_1,h_2,\ldots,h_n)\geq 0\quad \Longleftrightarrow
\quad
\varphi_k(t)\geq 0, t\in[-1,1], k=1,2,\ldots,n-1,
\]
where
\[
\varphi_k(t):=\left.f(h_1,h_2,\ldots,h_n)\right|_{h_1=h_2=\ldots=h_k=t,\,h_{k+1}=\ldots=h_n=1}
\]
for $k=1,2,\ldots,n-1$.
\end{thm}
In order to prove \eqref{buden2}, by Theorem \ref{propTi},
we only need to check $\varphi_k(t)\geq 0$ for all
$k=1,2,\ldots,n-1$. Indeed,
\[
\varphi_1(t)=(t-1)^2 t^2\geq 0,\quad\quad\varphi_{n-1}(t)=(t-1)^2\geq 0
\]
and
\[
\varphi_{k}(t)=\frac{k (n-k) (t-1)^2 }{(n-1)(n-2)^2}(a_1t^2+b_1t+c_1),\quad k=2,\ldots,n-2
\]
where
\begin{equation*}
\begin{aligned}
a_1:=&\left[(n-1)k-1\right](n-1-k)>0,\\
b_1:=& -2(k-1)(n-1-k)(n-1),\\
c_1:=& (k-1)\left[(n-k)(n-1)-1\right].
\end{aligned}
\end{equation*}
Since $a_1>0$ and
\[
b_1^2-4a_1c_1=-4(k-1)(n-1-k)n(n-2)^2<0,
\]
we conclude that
\[
a_1t^2+b_1t+c_1>0  \quad\mathrm{and}\quad
\varphi_{k}(t)\geq 0
\]
for $k=2,\ldots,n-2$. Therefore \eqref{buden2} follows.
\end{proof}

Using Lemmas \ref{tech1} and \ref{tech2},
we now finish the proof of the nonnegativity of $Q$.
\begin{theorem}\label{thm}
Let $\alpha_0$, $\alpha_k$, $s_k$ and $\beta$ be defined by \eqref{def1},
\eqref{def2} and \eqref{def3}. Then
\[
Q=\gamma^T(s_2I_n+\alpha_0\beta^T+\beta\alpha^T_0)\gamma\geq 0.
\]
\end{theorem}

\begin{proof}[Proof of Theorem \ref{thm}]
Let $B:=s_2I_n+\alpha_0\beta^T+\beta\alpha^T_0$. We only to show that
all eigenvalues of real matrix $B$ are real, nonnegative, and their
number is $n$. In fact, if $\lambda\neq s_2$, by Lemma
\ref{tech1}, we compute that
\begin{equation*}
\begin{aligned}
\det(\lambda I_n-B)&=\det((\lambda-s_2)I_n-\alpha_0\beta^T-\beta\alpha_0^T)\\
&=\det((\lambda-s_2)I_n)\cdot\det\left(I_2+\frac{1}{(\lambda-s_2)}(\beta,\alpha_0)^T(-\alpha_0,-\beta)\right)\\
&=(\lambda-s_2)^{n-2}\det\left((\lambda-s_2)I_2-(\beta,\alpha_0)^T(\alpha_0,\beta)\right)\\
&=(\lambda-s_2)^{n-2}(\lambda^2-p\lambda+q),
\end{aligned}
\end{equation*}
where
\[
p:=2(s_2+\alpha_0^T\beta)\quad\mathrm{and}
\quad
q:=(s_2+\alpha_0^T\beta)^2-n\beta^T\beta.
\]
In other words, when
$\lambda\neq s_2$, we have the identity
\begin{equation}\label{eigenequ}
\det(\lambda I_n-B)=(\lambda-s_2)^{n-2}(\lambda^2-p\lambda+q).
\end{equation}
Notice that two hand sides of the above identity are continuous
with respect to the parameter $\lambda$. Thus if $\lambda\to s_2$,
we know that \eqref{eigenequ} also holds for $\lambda=s_2$.
Therefore, \eqref{eigenequ} in fact holds \emph{without} any condition.
Now from \eqref{eigenequ}, we easily conclude that $s_2$ is
an nonnegative eigenvalue  of multiplicity $n-2$ in $B$.

On the other hand, by Lemma \ref{tech2}, we know that $p,\,q\geq 0$.
Meanwhile,
\[
p^2-4q=n\beta^T\beta\geq 0.
\]
Hence, equation $\lambda^2-p\lambda+q=0$ has two real solutions,
which implies that $B$ has another two nonnegative real eigenvalues.

In summary, we prove that the number of
nonnegative real eigenvalues of matrix $B$ is $n$. Therefore
$Q=\gamma^TB\gamma\geq 0$.
\end{proof}

Combining Theorem \ref{thm} and \eqref{PQ} implies that
\begin{corollary}\label{cor}
If $n\geq 4$, for any Riemannian metric $g$ and symmetric $2$-tensor $h$, we have
\[
|h|^2|Rc|^2-\frac{2R}{n-2}(HR_{ij}h_{ij}-R_{jk}h_{ji}h_{ik})+\frac{R^2(H^2-|h|^2)}{(n-1)(n-2)}\geq 0.
\]
\end{corollary}
\begin{remark}
We would like to point out that our proof method is also
suitable for the case $n=3$, which has been proved by
G. Anderson and B. Chow \cite{[AC]}.
\end{remark}

\bibliographystyle{amsplain}

\begin{thebibliography}{30}
\bibitem{[Per]} G. Perelman, The entropy formula for the Ricci flow and
its geometric applications, arxiv:0211159.
\end{thebibliography}

\end{document}